\documentclass[12pt]{amsart}
\usepackage{amsmath,amscd,amssymb,amsfonts}
\newcommand{\ver}{Feb. 2, 2007, v.9}
\newcommand{\ssbull}{\raise.2ex\hbox{${\scriptscriptstyle\bullet}$}}
\newcommand{\msum}{\hbox{$\sum$}}

\newcommand{\moplus}{\hbox{$\bigoplus$}}
\newcommand{\mtims}{\hbox{$\times$}}
\newcommand{\bC}{{\mathbb C}}
\newcommand{\bN}{{\mathbb N}}
\newcommand{\bP}{{\mathbb P}}
\newcommand{\bQ}{{\mathbb Q}}
\newcommand{\bZ}{{\mathbb Z}}
\newcommand{\bw}{{\mathbf w}}
\newcommand{\cA}{{\mathcal A}}
\newcommand{\cB}{{\mathcal B}}
\newcommand{\cD}{{\mathcal D}}
\newcommand{\cG}{{\mathcal G}}
\newcommand{\cH}{{\mathcal H}}
\newcommand{\cI}{{\mathcal I}}
\newcommand{\cJ}{{\mathcal J}}
\newcommand{\cL}{{\mathcal L}}
\newcommand{\cO}{{\mathcal O}}
\newcommand{\cS}{{\mathcal S}}
\newcommand{\hM}{\widehat{M}}
\newcommand{\hN}{\widehat{N}}

\newcommand{\oom}{\overline{\omega}}
\newcommand{\tH}{\widetilde{H}}
\newcommand{\tP}{\widetilde{P}}
\newcommand{\tU}{\widetilde{U}}
\newcommand{\tY}{\widetilde{Y}}
\newcommand{\tcG}{\widetilde{\mathcal G}}
\newcommand{\nnc}{\text{{\rm nnc}}}
\newcommand{\red}{\text{{\rm red}}}
\newcommand{\sm}{\text{{\rm sm}}}
\newcommand{\Gr}{\text{{\rm Gr}}}

\newcommand{\Sp}{\text{{\rm Sp}}}
\newcommand{\DR}{\text{{\rm DR}}}
\newcommand{\Sing}{\text{{\rm Sing}}\,}

\begin{document}
\title[Multiplier ideals, $b$-function, and spectrum]
{Multiplier ideals, $b$-function, and spectrum\\
of a hypersurface singularity}
\author{Morihiko Saito}
\address{RIMS Kyoto University, Kyoto 606-8502 Japan}
\email{msaito@kurims.kyoto-u.ac.jp}
\date{\ver}
\begin{abstract}
We prove that certain roots of the Bernstein-Sato polynomial (i.e.
$b$-function) are jumping coefficients up to a sign,
showing a partial converse of a theorem of L.~Ein, R.~Lazarsfeld,
K.E.~Smith, and D.~Varolin.
We also prove that certain roots are determined by a filtration on
the Milnor cohomology, generalizing a theorem of B.~Malgrange in the
isolated singularity case.
This implies a certain relation with the spectrum which is determined
by the Hodge filtration, because the above filtration is related to
the pole order filtration.
For multiplier ideals we prove an explicit formula in the case of
locally conical divisors along a stratification, generalizing a
formula of Musta\c{t}\v{a} in the case of hyperplane arrangements.
We also give another proof of a formula of U.~Walther on the
$b$-function of a generic hyperplane arrangement, including
the multiplicity of $-1$.
\end{abstract}
\maketitle

\centerline{\it To Joseph Steenbrink}

\bigskip
\centerline{\bf Introduction}

\bigskip\noindent
Let
$ X $ be a complex manifold, and
$ D $ be an effective divisor on it.
For a positive rational number
$ \alpha $, the multiplier ideal
$ \cJ(X,\alpha D) $ is a coherent ideal of the structure sheaf
$ \cO_{X} $ defined by the local integrability of
$ |g|^{2}/|f|^{2\alpha} $ for
$ g \in \cO_{X} $,
where
$ f $ is a holomorphic function defining
$ D $ locally, see [12], [21], [27].
This gives a decreasing filtration on
$ \cO_{X} $,
and essentially coincides with the filtration induced by the
$ V $-filtration of M.~Kashiwara [18] and
B.~Malgrange [25] along
$ D $ indexed by
$ \bQ $, see [5].
It is also related to the spectrum
$ \Sp(f,x) $, see [4], [5].

The numbers
$ \alpha $ at which the
$ \cJ(X,\alpha D) $ jump are called the jumping coefficients of
$ D $.
It is shown by L.~Ein, R.~Lazarsfeld, K.E.~Smith, and D.~Varolin
(see [12]) that any jumping coefficients which are less than
$ 1 $ are roots of the Bernstein-Sato polynomial
$ b_{f}(s) $ (i.e. the
$ b $-function) up to a sign.
It is well known that the minimal jumping coefficient
$ \alpha_{f} $ coincides with the minimal root of
$ b_{f}(-s) $, see [19].
For
$ x \in D $, we define
$ b_{f,x}(s) $,
$ \alpha_{f,x} $ by replacing
$ X $ with a sufficiently small neighborhood of
$ x $.
For
$ \alpha > 0 $ with
$ 0 < \varepsilon \ll 1 $, the graded pieces are defined by
$$
\cG(X,\alpha D) = \cJ(X,(\alpha-\varepsilon) D)/\cJ(X,\alpha D)\,
(= \Gr_{V}^{\alpha}\cO_{X}).
$$
We say that
$ \alpha $ is a local jumping coefficient of
$ D $ at
$ x $ if
$ \cG(X,\alpha D)_{x} \ne 0 $.
We have a partial converse to the theorem of [12] as follows
(see 3.3):

\medskip\noindent
{\bf Theorem~1.}
{\it Let
$ \alpha $ be a root of
$ b_{f,x}(-s) $ contained in
$ (0,1) $.
Assume

\noindent
$ (a) \,\, \xi f = f $ for a holomorphic vector field
$ \xi $.

\noindent
$ (b) \,\, \alpha < \alpha_{f,y} $ for any
$ y\ne x $ sufficiently near
$ x $.

\noindent
Then
$ \alpha $ is a local jumping coefficient of
$ D $ at
$ x $.
}

\medskip
Theorem~1 does not hold if either of the two conditions is not
satisfied, see 3.4.
Condition
$ (b) $ is satisfied if
$ \exp(-2\pi i \beta) $ is not an eigenvalue of the Milnor
monodromy of
$ f $ at
$ y \ne x $ for any
$ \beta \in [\alpha_{f,x},\alpha] $.
By definition,
$ \cJ(X,(\alpha+1)D) = f\cJ(X,\alpha D) $ for
$ \alpha > 0 $, and the jumping coefficients have a periodicity so
that
$ \alpha > 0 $ is a jumping coefficient if and only if
$ \alpha + 1 $ is.
However, the roots of
$ b_{f}(-s) $ do not have such a periodicity and we have to restrict
to
$ (\alpha_{f,x},1) $.

As for the relation with the spectrum,
N.~Budur [4] proved that, if
$ \alpha \in (0,1) $ and
$ \cG(X,\alpha D) $ is supported on a point
$ x $ of
$ D $ (e.g. if condition
$ (b) $ of Theorem~1 is satisfied),
then the coefficient
$ m_{\alpha} $ of the spectrum
$ \Sp(f,x) = \sum_{\beta} m_{\beta}t^{\beta} $ is given by
$$
m_{\alpha} = \dim \cG(X,\alpha D)_{x}.
\leqno(0.1)
$$
Indeed, under the above hypothesis,
$ \cG(X,\alpha D) \,(= \Gr_{V}^{\alpha}\cO_{X}) $ is identified
with the Hodge filtration
$ F^{n-1} $ on the
$ \lambda $-eigenspace of the Milnor monodromy
$ H^{n-1}(F_{x},\bC)_{\lambda} $ for
$ \lambda = \exp(-2\pi i\alpha) $ where
$ F_{x} $ denotes the Milnor fiber around
$ x $ and
$ n = \dim X $, see [5].
In the isolated singularity case, (0.1) is closely
related to [22], [24], [40], [42].
We have a generalization of a result of Malgrange [24] as follows
(see 4.5):

\medskip\noindent
{\bf Theorem~2.}
{\it There exists canonically a decreasing filtration
$ \tP $ on
$ H^{n-1}(F_{x},\bC)_{\lambda} $ stable by the Milnor monodromy and
containing the Hodge filtration
$ F $, and for any rational number
$ \alpha $ such that
$ \lambda = \exp(-2\pi i\alpha) $, we have the following\,{\rm :}

\noindent
$ (a) $
If
$ \Gr_{\tP}^{p}H^{n-1}(F_{x},\bC)_{\lambda} \ne 0 $ with
$ p = [n-\alpha] $, then
$ \alpha $ is a root of
$ b_{f,x}(-s) $.

\noindent
$ (b) $
If
$ \alpha + i $ is not a root of
$ b_{f,y}(-s) $ for any
$ y \ne x $ and any
$ i \in \bN $, then the converse of the assertion
$ (a) $ holds.

\noindent
$ (c) $
If
$ \lambda $ is not an eigenvalue of the Milnor monodromy at
$ y\ne x $, then the multiplicity of the root
$ \alpha $ coincides with
the degree of the minimal polynomial of the action of the monodromy on
$ \Gr_{\tP}^{p}H^{n-1}(F_{x},\bC)_{\lambda} $.
}

\medskip
Note that the spectrum [39] is defined by the same way as in
$ (a), (c) $ replacing
$ \tP $ with the Hodge filtration
$ F $ and the minimal polynomial with the characteristic polynomial,
see 3.5.
The filtration
$ \tP $ is defined by using the saturated Brieskorn lattices
$ \tcG_f^{(-i)} $ (see (4.1.6)), and contains the Hodge filtration
$ F $, see Proposition~4.4.
Replacing
$ \tcG_f^{(0)} $ with the Brieskorn lattices
$ \cG_f^{(-i)} $, we have the filtration
$ P $ contained in
$ \tP $, see (4.1.6).
In the isolated singularity case,
$ P $ coincides with the Hodge filtration
$ F $, see [36], [42].
In the quasi-homogeneous isolated singularity case, this also follows
from [24], [37]
(where the Milnor cohomology is identified with the Jacobian ring, and
the Hodge filtration is described by using the weighted degree of
monomials).
If
$ f $ is a homogeneous polynomial in general, then
$ \tP $ coincides with
$ P $ and with the pole order filtration defined by using a local system
on an open subvariety of
$ \bP^{n-1} $ calculating
$ H^{n-1}(F_{x},\bC)_{\lambda} $, see Proposition~4.9.

In general it is not easy to calculate
$ \cJ(X,\alpha D) $ explicitly except for some special cases, see
[15], [16], [26], etc.
In this paper, we give an explicit formula for
$ \cJ(X,\alpha D) $ in the case
$ D $ is a locally conical divisor along a stratification, i.e.
$ D $ is locally defined by a weighted homogeneous function
with nonnegative weights and the zero weight part,
which is the limit of the (local)
$ \bC^{*} $-action, is given by the stratum passing through the point
(e.g.
$ D $ is an affine cone of a divisor on
$ \bP^{2} $ which is defined locally in classical topology
by a weighted homogeneous polynomial),
see 1.2 for details.
We have a (shifted) decreasing filtration
$ \{G_{x}^{\alpha}\}_{\alpha\in\bQ} $ on
$ \cO_{X,x} $ associated to the weights at each
$ x \in D $, see 1.3.
Let
$ D^{\nnc} $ denote the smallest closed analytic subset of
$ D $ such that
$ D $ is a divisor with normal crossings outside
$ D^{\nnc} $, and let
$ D_{\red}^{\sm} $ be the smooth part of the reduced variety
$ D_{\red} $.

\medskip\noindent
{\bf Theorem~3.}
{\it Let
$ X $ be a complex manifold, and
$ D $ be a locally conical divisor along a stratification.
Then a section
$ g $ of
$ \cO_{X} $ belongs to
$ \cJ(X,\alpha D) $ if and only if
$ g_{x} \in \cO_{X,x} $ belongs to
$ G_{x}^{>\alpha} $ for any
$ x \in D^{\nnc} \cup D_{\red}^{\sm} $.
}

\medskip
This generalizes a formula of Musta\c{t}\v{a} [26] for a
hyperplane arrangement with a reduced equation.
The condition for
$ x \in D_{\red}^{\sm} $ is equivalent to that the vanishing
order of
$ g $ along
$ D $ is strictly greater than
$ \alpha m_{x} - 1 $ where
$ m_{x} $ is the multiplicity of
$ D $ at
$ x $.
A similar formula has been known for a function with
nondegenerate Newton boundary,
see [15], [16], [21] (and also 2.5 below).
By induction on stratum, Theorem~3 is reduced to Theorem~2.2 below
whose proof uses the above analytic definition
of multiplier ideal together with some commutative algebra,
see 2.1--3.

For a divisor
$ D $ on a complex manifold, let
$ \alpha_{D} = \min\{\alpha_{f,x} : x \in D\} $ where
$ f $ is a holomorphic function defining
$ D $ on a neighborhood of
$ x $.
As a corollary of Theorem~2.2 we can deduce

\medskip\noindent
{\bf Proposition~1.}
{\it Assume
$ X = \bC^{n} $ and
$ D $ is the affine cone of a divisor
$ Z $ of degree
$ d $ on
$ \bP^{n-1} $.
Let
$ \cI_{0} $ be the ideal sheaf of
$ \{0\} \subset \bC^{n} $.
Then we have for
$ \alpha < \alpha_{Z} $
$$
\cJ(X,\alpha D) = \cI_{0}^{k} \quad\text{with}\quad
k = [d\alpha] - n + 1.
$$
In particular,
$ j/d $ is a local jumping coefficient of
$ D $ at
$ 0 $ if
$ n \le j < d\alpha_{Z} $.
}

\medskip
Note that
$ \alpha_{Z} \le 1 $, and the equality holds if
$ Z $ is a reduced divisor with normal crossings, e.g. if
$ D $ is a generic hyperplane arrangement, see also [26].
Since
$ \dim \cI_{0}^{k}/\cI_{0}^{k+1} = \binom{n+k-1}{n-1} $,
we get from (0.1) and Proposition~1 the following

\medskip\noindent
{\bf Corollary~1.}
{\it With the notation and the assumption of Proposition {\rm 1},
assume
$ Z $ is a reduced divisor with normal crossings on
$ \bP^{n-1} $.
Then the coefficients
$ m_{\alpha} $ and
$ m_{n-\alpha} $ of the spectrum
$ \Sp(f,0) $ are
$ \binom{j-1}{n-1} $ for
$ \alpha = j/d < 1 $.
}

\medskip
The assertion for
$ m_{n-\alpha} $ is reduced to that for
$ m_{\alpha} $ by the symmetry of the Hodge numbers
for the nonunipotent monodromy part of the vanishing cohomology
(which is identified with that of the nearby cycle sheaf in this
case).
Note that the formula is the same as in the case of a homogeneous
polynomial with an isolated singularity, and can also be deduced
from the calculation of the Hodge filtration in 4.8.

In the case of a generic central hyperplane arrangement
with a reduced equation
$ f $, the
$ b $-function is determined by U.~Walther [43] (except for the
multiplicity of the root
$ - 1 $).
Using Theorem~2 together with [3], [14], we first see that
the roots of
$ b_{f}(-s) $ are strictly smaller than
$ 2 $, see Proposition~5.2.
Then we can give another proof of his formula
together with the multiplicity of
$ -1 $, using the relation between the
$ b $-function and the
$ V $-filtration in [18], [25] together with Corollary~1,
see 5.4.
Note that for any arbitrary hyperplane arrangement,
$ -1 $ is the only integral root of
$ b_{f}(s) $ (see [43]), and we can show that its multiplicity is
$ -n $ if the arrangement is not the pull-back of an arrangement in
a strictly lower dimensional vector space, see Proposition~5.3.
More details will be given in a forthcoming paper on the
$ b $-functions of hyperplane arrangements.

Walther's formula shows that, without restricting to the interval
$ (0,1) $, there is no relation between the spectrum and the roots of
$ b_{f}(-s) $ (contrary to the case of a homogeneous polynomial with
an isolated singularity).
This comes from the difference between the Hodge and pole order
filtrations on the Milnor cohomology in Theorem~2, see Proposition~4.9.

As for the spectrum of generic central hyperplane arrangement, the
$ m_{\alpha} $ for
$ \alpha \in \bZ $ are easily calculated, see (5.6.1).
Combined with Corollary~1, this gives the spectrum of a generic
central hyperplane arrangement for
$ n = 3 $, because the Euler characteristic is calculated in
[6], [28].
It is possible in principle to calculate the spectrum for a general
$ n $, using [10].

I would like to thank Dimca for useful discussions related to this paper,
and the referee for good comments.

In Sect.~1, we introduce locally conical divisors along a
stratification.
In Sect.~2, we prove Theorem~2.2 which implies Theorem~3.
In Sect.~3, we explain the relation with
$ b $-function and spectrum, and prove Theorem~1.
In Sect.~4, we explain the relation with Brieskorn modules and
Gauss-Manin systems, and prove Theorem~2.
In Sect.~5, we treat the case of a generic central hyperplane
arrangement, and give another proof of Walther's theorem as an
application of Theorem~2.

\bigskip\bigskip
\centerline{{\bf 1. Locally conical divisors along a stratification}}

\bigskip\noindent
{\bf 1.1.~Conical divisors.}
Let
$ Y $ be a complex manifold, and
$ (x_{1}, \cdots, x_{r}) $ be the coordinate system of
$ \bC^{r} $.
Let
$ \bw = (w_{1}, \cdots, w_{r}) \in \bQ_{>0}^{r} $ (i.e.
$ w_{i} $ are positive rational numbers).
We say that an effective divisor
$ D $ on
$ X := Y\mtims \bC^{r} $ is a {\it conical divisor} along
$ Y\mtims\{0\} $ with
positive weight
$ \bw = (w_{1}, \cdots, w_{r}) $ if
$ D $ is locally defined by a relatively weighted homogeneous function
$ f $ with positive weight
$ \bw $ (i.e.
$ f $ is a linear combination of
$ x_{1}^{{\nu}_{1}} \cdots x_{r}^{{\nu}_{n}} $ with
$ \cO_{Y} $-coefficients such that
$ (\nu_{1}, \cdots, \nu_{r}) \in \bN^{r} $ satisfies
$ \sum_{i} w_{i}\nu_{i} = 1 $).
For a positive real number
$ \lambda $, we define
$ \phi_{\lambda} : \bC^{r} \to \bC^{r} $ by
$$
\phi_{\lambda}(x_{1}, \cdots, x_{r}) = (\lambda^{w_{1}}
x_{1}, \cdots, \lambda^{w_{r}}x_{r}),
\leqno(1.1.1)
$$
and
$ id \mtims \phi_{\lambda} : Y\mtims \bC^{r} \to Y\mtims
\bC^{r} $ will be denoted also by
$ \phi_{\lambda} $.
Then
$ {\phi}_{\lambda}^{*}f = \lambda f $.

\medskip
\noindent
{\bf 1.2.~Locally conical divisors along a stratification.}
We say that a divisor
$ D $ on a complex manifold
$ X $ is a {\it locally conical divisor} along a smooth
submanifold
$ Z $ if for each
$ z \in Z $,
there exist a complex manifold
$ Y $, a conical divisor
$ D' $ on
$ Y\mtims \bC^{r} $ with positive weight
$ \bw $ along
$ Y\mtims \{0\} $, an open subset
$ U' $ of
$ Y\mtims \bC^{r} $ and an open neighborhoods
$ U $ of
$ z $ in
$ X $ together with an isomorphism
$ U \simeq U' $ inducing isomorphisms
$ U\cap D \simeq U'\cap D' $,
$ U\cap Z \simeq U'\cap (Y\mtims \{0\}) $ (in particular,
$ z $ corresponds to a point of
$ Y\mtims \{0\} $).
Note that the weight
$ \bw $ is not necessarily unique in general.

Let
$ D $ be an effective divisor on a complex manifold
$ X $.
Let
$ D^{\nnc} $ be the smallest closed analytic subset of
$ D $ such that
$ D \setminus D^{\nnc} $ is a divisor with normal crossings
on
$ X \setminus D^{\nnc} $.
We say that
$ D $ is a {\it locally conical divisor along a stratification}
$ \{S_{i}\} $ of
$ D^{\nnc} $, if
$ D $ is a locally conical divisor along
$ S_{i} $ for each
$ i $.

\medskip
\noindent
{\bf 1.3.~Shifted $ \bw $-filtrations along strata.}
With the notation of 1.1, let
$ x = (y,0) \in Y\mtims \{0\} $ and
$ g \in \cO_{X,x} $.
We have the expansion
$$
g = \msum_{\beta} \,g_{\beta}
\leqno(1.3.1)
$$
such that
$ g_{\beta} $ is a linear combination of
$ x_{1}^{{\nu}_{1}} \cdots x_{r}^{{\nu}_{n}} $ with
$ \cO_{Y} $-coefficients satisfying
$$
\msum_{i} w_{i}(\nu_{i}+1) = \beta.
\leqno(1.3.2)
$$
We define a decreasing filtration
$ G $ of ideals of
$ \cO_{X,x} $ such that
$ G^{\alpha} $ is generated by
$ x_{1}^{{\nu}_{1}} \cdots x_{r}^{{\nu}_{n}} $ with
$ \sum_{i} w_{i}(\nu_{i}+1) \ge \alpha $
(i.e.
$ g \in G^{\alpha} $ if and only if
$ g_{\beta} $ vanishes for
$ \beta < \alpha $).
This is called the shifted
$ \bw $-filtration.

If
$ D $ is a locally conical divisor along
$ Z $ or a stratification
$ \{S_{i}\} $ as in (1.2-3), we have the shifted
$ \bw $-filtration
$ G_{x} $ on
$ \cO_{X,x} $ for each
$ x \in Z $ or
$ D^{\nnc} $.
This is not necessarily unique in general.

If
$ x \in D_{\red}^{\sm} \,(:= D \setminus \Sing D_{\red}) $,
let
$ h $ be a holomorphic function defining
$ D_{\red} $ on a neighborhood of
$ x $, and
$ m_{x} $ be the multiplicity of
$ D $ at
$ x $.
Then for
$ \alpha > 0, $ we have
$ G_{x}^{\alpha} = h^{i-1}\cO_{X,x} $ where
$ i $ is the minimal integer such that
$ i \ge m_{x}\alpha $.

Let
$ G_{x}^{>\alpha} = \bigcup_{\beta >\alpha} G_{x}^{\beta} $
in general.

\medskip\noindent
{\bf 1.4.~Remarks.}
(i) If
$ f = \sum_{i}u_{i}x_{i}^{a_{i}} $ with
$ u_{i}(0) \ne 0 $ and
$ a_{i} \in \bZ_{>0} $, then
$ D := f^{-1}(0) $ is locally conical along the origin,
because
$ f = \sum_{i}y_{i}^{a_{i}} $ with
$ y_{i} = u_{i}^{-1/a_{i}}x_{i} $.
This implies, for example, that a divisor
$ D $ on
$ \bC^{3} $ is locally conical along a stratification,
if it is defined by
$ f = x^{a}y^{b} + y^{a}z^{b} + z^{a}x^{b} $ with
$ a, b > 1 $.

\medskip
(ii) It is possible that the moduli of singularity really
vary along a stratum, e.g.
$ f = (x^{3}+y^{3}+z^{3})u + xyzv $.

\bigskip\bigskip
\centerline{{\bf 2. Calculation of multiplier ideals}}

\bigskip\noindent
The following is the key to the proof of Theorem~3.

\medskip\noindent
{\bf 2.1.~Proposition.}
{\it With the notation and the assumption of {\rm 1.1}, let
$ X = Y\mtims \bC^{r}, X' = X \setminus Y\mtims \{0\} $
with the inclusion
$ j : X' \to X $.
Put
$ D' = D \cap X' $.
If
$ g \in \cO_{X,x} $ belongs to
$ j_{*}\cJ(X',\alpha D') $,
then each
$ g_{\beta} $ in {\rm (1.3.1)} belongs to
$ j_{*}\cJ(X',\alpha D') $.
}

\medskip\noindent
{\it Proof.}
Since
$ \cJ(X',\alpha D') $ is extended to a coherent sheaf on
$ X $,
$ j_{*}\cJ(X',\alpha D') \cap \cO_{X} $ is coherent.
(Indeed, the assertion is reduced to the case where the
complement of the image of
$ j $ is a divisor, using a Cech complex.
Then any section of
$ j_{*}\cJ(X',\alpha D') \cap \cO_{X} $ defines a section of
$ \cO_{X}/\cJ(X,\alpha D) $, which is supported on the divisor,
and hence is annihilated by a sufficiently high power of a
function defining the divisor.
So the intersection with
$ j_{*}\cJ(X',\alpha D') $ can be replaced by the one with the
algebraic localization of a coherent extension of
$ \cJ(X',\alpha D') $ which is quasi-coherent.)

For
$ x = (y,0) \in Y \mtims \{0\} $, let
$$
M = (j_{*}\cJ(X',\alpha D') \cap \cO_{X})_{x},\quad
N = \cO_{X,x}.
$$
We will denote
$ G_{x} $ by
$ G $ in this subsection to simplify the notation.
For
$ \beta \in \bQ $, let
$$
G^{\beta}\hM = \varprojlim{}_{\gamma} G^{\beta}M/G^{\gamma}M,\quad
G^{\beta}\hN = \varprojlim{}_{\gamma} G^{\beta}N/G^{\gamma}N,
$$
and
$ \hM, \hN $ be their inductive limit for
$ \beta $ respectively.
By the Mittag-Leffler condition, we have the injectivity of
$ G^{\beta}\hM \to G^{\beta'}\hM $ for
$ \beta > \beta' $ so that we get the filtration
$ G $ of
$ \hM $ (similarly for
$ \hN $), see also [32].
By the Artin-Rees lemma,
$ G^{\beta}\hM, G^{\beta}\hN $ coincide with the
$ I $-adic completion of
$ G^{\beta}M, G^{\beta}N $ by the ideal
$ I $ of
$ Y\mtims \{0\} $, because the filtration
$ G $ is induced by
$ G $ on
$ N $ which is essentially equivalent to the
$ I $-adic filtration (i.e. there are positive rational numbers
$ \alpha, \beta $ such that
$ G^{i\alpha} \subset I^{i} \subset G^{i\beta} $ for
$ i \gg 0 $).

For
$ \lambda\in \bC^{*} $, we see that
$ M $,
$ G^{\beta}N $, and hence
$ G^{\beta}M $ are stable by the action of
$ \phi_{\lambda}^{*} $.
So the filtration
$ G $ on
$ \hM $ splits canonically (because
$ G $ on
$ G^{\beta}/G^{\gamma} $ does).
Thus, for
$ g \in G^{\beta}M $,
we have
$ g = g_{\beta} + g' $ where
$ g_{\beta} $ is as in (1.3.1) and
$ g' \in G^{>\beta}\hM \cap N $
because
$ g_{\beta} \in N $.
So the assertion is reduced to
$$
G^{>\beta}M = G^{>\beta}\hM \cap N,
\leqno(2.1.1)
$$
because this implies
$ g_{\beta} \in M $ so that we can proceed by induction on
$ \beta $ replacing
$ g $ with
$ g' $.

For the proof of (2.1.1), consider the commutative diagram
$$
\CD
0 @>>> G^{>\beta}M @>>> N @>>> N/G^{>\beta}M
@>>> 0
\\
@. @VVV @VVV @VVV
\\
0 @>>> (G^{>\beta}M)^{\wedge} @>>> \hN @>>>
(N/G^{>\beta}M)^{\wedge} @>>> 0
\endCD
$$
where the bottom row is the
$ I $-adic completion of the top row.
By the above argument, we have
$ G^{>\beta}\hM = (G^{>\beta}M)^{\wedge} $,
and the vertical morphisms are injective by Krull's
intersection theorem.
So (2.1.1) follows.
This completes the proof of Proposition~2.1.

\medskip\noindent
{\bf 2.2.~Theorem.}
{\it With the above notation and assumption, we have
$$
\cJ(X,\alpha D)_{x} = (j_{*}\cJ(X',\alpha D'))_{x} \cap
G^{>\alpha}\cO_{X,x}\quad \text{for}\,\,\,x \in Y\mtims \{0\},
\leqno(2.2.1)
$$
where
$ G^{>\alpha} $ is as in {\rm 1.3}.
}

\medskip\noindent
{\it Proof.}
We first show that
$ g \in \cJ(X,\alpha D) $ if
$ g \in (j_{*}\cJ(X',\alpha D'))_{x} \cap G^{>\alpha}\cO_{X,x} $.
We have the expansion
$ g = \sum_{\beta} g_{\beta} $ as in (1.3.1).
By Proposition~2.1 we may assume
$ g = g_{\beta} $ with
$ \beta > \alpha $, because
$ (j_{*}\cJ(X',\alpha D') \cap \cO_{X})/\cJ(X,\alpha D) $ is annihilated
by a sufficiently high power of the ideal of
$ Y\mtims \{0\} $ in the notation of Proposition~2.1
so that the assertion is clear if
$ g \in G_{x}^{\beta} $ for
$ \beta $ sufficiently large.
We have a relatively compact open subset
$ U $ of
$ X' $ together with
$ \lambda \in (0,1) $ such that
$ \bigcup_{j\ge 0}\phi_{\lambda^{j}}U $ contains
$ U_{x} \setminus Y\mtims\{0\} $ where
$ U_{x} $ is a sufficiently small open neighborhood of
$ x $ in
$ X $ on which
$ g $ is defined.
(For example, use a function defined by
$ \rho (x) = \sum_{i\le r} |x_{i}|^{1/w_{i}} $ so that
$ \rho (\phi_{\lambda}x) = \lambda \rho (x) $.)

Let
$ \omega = dx_{1}\wedge \cdots\wedge dx_{r}\wedge dy_{1}
\wedge \cdots\wedge dy_{n-r} $,
where
$ (y_{1}, \cdots, y_{n-r}) $ is a local coordinate system
of
$ Y $.
Then
$ {\phi}_{\lambda^{j}}^{*}(g_{\beta}\omega) =
\lambda^{j\beta}g_{\beta}\omega $,
and
$$
\aligned
\msum_{j\ge 0}\int_{U\setminus D}{\phi}_{{\lambda}^{j}}^{*}
((|g_{\beta}|^{2}/|f|^{2\alpha})\omega \wedge \oom)
&= \msum_{j\ge 0}\lambda^{2(\beta -\alpha)j}\int_{U\setminus D}
(|g_{\beta}|^{2}/|f|^{2\alpha})\omega \wedge \oom
\\
&= (1-\lambda^{2(\beta -\alpha)})^{-1}\int_{U\setminus D}
(|g_{\beta}|^{2}/|f|^{2\alpha})\omega \wedge \oom.
\endaligned
$$
So the assertion follows.

Similarly, we see that
$ g \notin \cJ(X,\alpha D) $ if
$ g \in (j_{*}\cJ(X',\alpha D'))_{x} $ and
$ g \notin G_{x}^{>\alpha} $.
Here we may assume
$ g \in G_{x}^{\alpha} $ by replacing
$ \alpha $ with a smaller number if necessary.
Then we may assume further that
$ g = g_{\alpha} $ using the above argument.
So the assertion follows by
considering a sufficiently small open subset
$ U $ of
$ X' $ such that the
$ \phi_{\lambda^{j}}U \,(j \in \bN) $ are disjoint.
This completes the proof of Theorem~2.2.

\medskip\noindent
{\bf 2.3.~Proof of Theorem~3.}
The assertion is well-known outside
$ D^{\nnc} $, i.e. if
$ D $ is a divisor with normal crossings, see e.g. [4], [5].
We proceed by induction on stratum.
Since the assertion is local, we may assume
$ X = Y\mtims \bC^{r} $ with
$ Y = S_{i} $ and
$ D $ is defined by a relatively homogeneous function
$ f $ with positive weight
$ \bw = (w_{1}, \cdots, w_{r}) $ as in 1.1.
Then the assertion follows from Theorem~2.2 by
induction on stratum.

\medskip\noindent
{\bf 2.5.~Nondegenerate Newton boundary case.}
Assume
$ f \in \cO_{X,x} $ has a nondegenerate Newton boundary
([20], [41]).
Then we have a formula similar to Theorem~3
by [15], [16], [21].
There is a shifted Newton filtration
$ G $ on
$ \cO_{X,0} $ such that
$ G^{\alpha}\cO_{X,0} $ is generated over
$ \cO_{X,0} $ by the monomials
$ x_{1}^{\nu_{1}}\cdots x_{n}^{\nu_{n}} $ satisfying
$$
\msum_{i} w_{\sigma,i}(\nu_{i}+1) \ge \alpha
\leqno(2.5.1)
$$
for any
$ (n-1) $-dimensional faces
$ \sigma $ of the Newton polyhedron, where the
$ w_{\sigma,i} $ are positive rational numbers such that
$ \sigma $ is contained in the hyperplane defined by
$ \sum_{i} w_{\sigma,i}\nu_{i} = 1 $.
Then
$$
\cJ(X,\alpha D)_{x} = G^{>\alpha}\cO_{X,0}
\quad\text{for}\,\, \alpha < 1.
\leqno(2.5.2)
$$
This is proved in loc.~cit. in the polynomial case.
The proof in the analytic case should be essentially
same.
(It would also be possible to use an argument similar to the
proof of Theorem~2.2
together with the torus embedding constructed in [41],
because the nondegeneracy corresponds to the condition that
the restriction of the proper transform of the hypersurface to
each stratum, which is isomorphic to a torus, is nonsingular.)
In the isolated singularity case with nondegenerate Newton
boundary, this is related to [32] using 3.2 and Proposition~4.7 below.

\bigskip\bigskip
\centerline{{\bf 3. Relation with
$ b $-function and spectrum}}

\bigskip\noindent
{\bf 3.1.~$b$-Function.}
Let
$ X $ be a complex manifold of dimension
$ n $, and
$ f $ be a non-constant holomorphic function on
$ X $.
Let
$$
M = \cD_{X}[s]f^{s}.
$$
It is identified with a
$ \cD_{X}[s] $-submodule of
$$
\cB_{f} := \cO_{X}\otimes_{\bC}\bC[\partial_{t}],
$$
generated by
$ 1\otimes 1 $ (which is identified with
$ f^{s} $), where
$ s = -\partial_{t}t $, see [17], [25].
Here
$ \cB_{f} $ is the direct image of
$ \cO_{X} $ by the graph embedding
$ i_{f} : X\to X\times \bC $ as a
$ \cD $-module, and the action of
$ \cD_{X\times\bC}$ on
$ \cB_{f} $ is defined by identifying
$ 1\otimes 1 $ with the delta function
$ \delta(t-f) $.
More precisely, for a vector field
$ \xi $ on
$ X $ and the coordinate
$ t $ of
$ \bC $, we have
$$
\aligned
\xi(g\otimes\partial_{t}^{j})
&=(\xi g)\otimes\partial_{t}^{j}-
(\xi f)g\otimes \partial_{t}^{j+1},
\\
t(g\otimes\partial_{t}^{j})
&=fg\otimes\partial_{t}^{j}-
jg\otimes \partial_{t}^{j-1},
\endaligned
\leqno(3.1.1)
$$
and the actions of
$ h \in \cO_{X} $ and
$ \partial_{t}^{i} $ are natural ones, see also [5].

The
$ b $-function
$ b_{f}(s) $ is the minimal polynomial of the action of
$ s $ on
$ M/tM $.
Since
$ M/tM $ is holonomic, the
$ b $-function exists if
$ X $ is (relatively) compact or
$ X, f $ are algebraic.
By M.~Kashiwara [18] and B.~Malgrange [25],
$ \cB_{f} $ has the filtration
$ V $ together with a canonical isomorphism of perverse sheaves
$$
\DR_{X}(\Gr_{V}^{\alpha}\cB_{f}) = \psi_{f,\lambda}\bC_{X}[n-1]
\quad\text{for}\,\,\,\alpha > 0,\,\lambda = \exp(-2\pi i\alpha)
\leqno(3.1.2)
$$
such that
$ \exp(-2\pi i\partial_{t}t) $ on the left-hand side
corresponds to the monodromy
$ T $ on the right-hand side.
Here
$ \DR_{X} $ denotes the de Rham functor (which induces an
equivalence of categories between regular holonomic
$ \cD $-modules and perverse sheaves) and
$ \psi_{f,\lambda}\bC_{X}[n-1] $ is the
$ \lambda $-eigenspace of the nearby cycle (perverse) sheaf
$ \psi_{f}\bC_{X}[n-1] $ for the semisimple part of the monodromy
$ T $, see [2], [9].

\medskip\noindent
{\bf 3.2.~Relation of the multiplier ideals with the
$ V $-filtration.}
By [5] we have
$$
\cJ(X,\alpha D) = V^{\alpha}\cO_{X} \quad\text{if
$ \alpha $ is not a jumping coefficient,}
\leqno(3.2.1)
$$
where the filtration
$ V $ on
$ \cO_{X} $ is induced by the
$ V $-filtration on
$ \cB_{f} \,(= \cO_{X}\otimes \bC[\partial_{t}]) $ in 3.1.
If
$ \alpha $ is a jumping coefficient (or actually, for any
$ \alpha $),
we have for
$ 0 < \varepsilon \ll 1 $
$$
\cJ(X,\alpha D) = V^{\alpha +\varepsilon}\cO_{X}, \quad
V^{\alpha}\cO_{X} = \cJ(X,(\alpha -\varepsilon)D).
\leqno(3.2.2)
$$

This implies another proof of a theorem of L.~Ein, R.~Lazarsfeld,
K.E.~Smith, and D.~Varolin (see [12]) that any jumping coefficients
which are less than
$ 1 $ are roots of the
$ b $-function up to a sign.

\medskip\noindent
{\bf 3.3.~Proof of Theorem~1.}
By 3.2 we can essentially replace
$ \cJ(X,\alpha D) $ with
$ V^{\alpha}\cO_{X} $.
By condition
$ (a) $, we have
$ \xi f = f $ so that
$ \xi f^{s} = sf^{s} $, and hence
$$
M := \cD_{X}[s]f^{s} = \cD_{X}f^{s} \subset \cB_{f}.
\leqno(3.3.1)
$$
By condition
$ (b) $,
$ M/V^{>\alpha}M $
is supported on
$ \{x\} $, and is generated over
$ \bC[\partial] := \bC[\partial_{1}, \cdots, \partial_{n}] $ by
$$
(\cO_{X}/V^{>\alpha}\cO_{X})\otimes 1 \subset
\cB_{f}/V^{>\alpha}\cB_{f},
\leqno(3.3.2)
$$
where
$ \partial_{j} = \partial /\partial x_{j} $.
Consider the filtered morphism induced by (3.3.2)
$$
(\cO_{X}/V^{>\alpha}\cO_{X},V)\otimes_{\bC}\bC[\partial]\to
(\cB_{f}/V^{>\alpha}\cB_{f},V).
\leqno(3.3.3)
$$
This is strictly injective,
i.e. it induces injective morphisms of the graded pieces.
Indeed, for
$ \beta \le \alpha $, the
$ \Gr_{V}^{\beta}\cO_{X,x} $ are finite dimensional vector spaces,
and are annihilated by the maximal ideal of
$ \cO_{X,x} $, see [4] and [33], 3.2.6.
Moreover,
$ V $ on
$ \cB_{f}/V^{>\alpha}\cB_{f} $ is a filtration as
$ \cD_{X} $-modules and the morphism of graded pieces induced by
(3.3.2) is injective by the definition of the induced filtration.
So we get the strictly injectivity of (3.3.3),
because any holonomic
$ \cD_{X} $-module supported on
$ \{x\} $ is isomorphic to a direct sum of
$ \bC[\partial] $, and is freely generated over
$ \bC[\partial] $ by its annihilator of the maximal ideal of
$ \cO_{X,x} $.

By the above argument, the image of (3.3.3) is
$ M/V^{>\alpha}M $, and the strict injectivity implies
that there is no
$ \beta \le \alpha $ such that
$ \Gr_{V}^{\beta} M \ne 0 $ but
$ \Gr_{V}^{\beta} \cO_{X} = 0 $.
So the assertion follows.

\medskip\noindent
{\bf 3.4.~Remarks.} (i)
The assertion of Theorem~1 does not hold if either of the
two conditions is not satisfied.
For example, consider
$ f = x^{5} + y^{4} + x^{3}y^{2} $ or
$ f = x^{5} + y^{4} + x^{3}y^{2}z $ where condition
$ (a) $ or
$ (b) $ is not satisfied respectively.
Here
$ \alpha = 4/5+3/4-1 = 11/20 > \alpha_{f} = 1/5 + 1/4 = 9/20 $.
Note that their jumping coefficients coincide with those for
$ f = x^{5} + y^{4} $ by [15], and
$ 11/20 $ is not a jumping coefficient.

\medskip
(ii) In the case of
$ f = (x^{2}-y^{2})(x^{2}-z^{2})(y^{2}-z^{2})z $, we see that
$ 5/7 $ is not a jumping coefficient by an argument in [26] (because
there is no hypersurface of degree
$ 2 $ on
$ \bP^{2} $ containing all the points of
$ Z $ whose multiplicity is
$ 3 $), but it is a root of
$ b_{f}(-s) $ as shown in 5.5 below.
In this case condition
$ (a) $ in Theorem~1 with positive weights is satisfied,
but condition
$ (b) $ is not.

\medskip
(iii) The first assertion of Theorem~1 trivially follows
from [25], if any rational number
$ \beta $ in
$ (\alpha_{f},\alpha) $ such that
$ \exp(-2\pi i\beta) $ is an eigenvalue of the Milnor
monodromy is a jumping coefficient.
This condition for any
$ \alpha \in (\alpha_{f},1) $ is satisfied by generic central
hyperplane arrangement, but not necessarily by nongeneric ones,
e.g. if
$ f $ is as in Remark~(ii) above.

\medskip\noindent
{\bf 3.5.~Spectrum.}
With the notation of 3.1, let
$ F_{x} $ denote the Milnor fiber around
$ x \in D := f^{-1}(0) $.
The spectrum
$ \Sp(f,x) = \msum_{\alpha \in \bQ} m_{\alpha}t^{\alpha} $
is defined by
$$
\aligned
&m_{\alpha} = \msum_{j} (-1)^{j-n+1} \dim \Gr_{F}^{p}
\tH^{j}(F_{x},\bC)_{\lambda}
\\
&\quad\text{with}\,\, p = [n - \alpha],\,\,
\lambda = \exp(-2\pi i\alpha),
\endaligned
$$
where
$ \tH^{j}(F_{x},\bC)_{\lambda} $ is the
$ \lambda $-eigenspace of the reduced cohomology for the
semisimple part
$ T_{s} $ of the Milnor monodromy
$ T $, and
$ F $ is the Hodge filtration, see [38], [39].

In this paper we define a mixed Hodge structure [8]
on the Milnor cohomology
$ H^{j}(F_{x},\bC) $ by using the pull-back of the nearby cycle sheaf
$ \psi_{f}\bC_{X}[n-1] $ by the inclusion
$ i_{x} : \{x\} \to X $ in the derived category of mixed Hodge modules.
This pull-back is defined by iterating the pull-back by
$ i_{j} : X^{j} \to X^{j-1} $, where
$ X^{j} = \{z_{i} = 0 : i \le j\} \subset X $ with
$ z_{1}, \dots, z_{n} $ local coordinates around
$ x $.
(Here we may assume that
$ X $ is a polydisk around
$ x $.)
The pull-back by
$ i_{j} $ is defined by using the mapping cone of
$ \partial_{j}:=\partial/\partial z_{j} : \Gr_{V_{j}}^{1} \to
\Gr_{V_{j}}^{0} $ where
$ V_{j} $ is the
$ V $-filtration of Kashiwara and Malgrange along
$ \{z_{j} = 0 \} $ and the Hodge filtration
$ F $ on
$ \Gr_{V_{j}}^{1} $ is shifted by
$ 1 $ so that
$ \partial_{j} $ preserves
$ F $.
(We can prove (0.1) using this, because
$ \cG(X,\alpha D) $ is annihilated by the maximal ideal under the
assumption of the (0.1).)

\medskip
The following lemma will be used in Proposition~5.3 to determine the
multiplicity of the root
$ -1 $ of the
$ b $-function of a hyperplane arrangement.

\medskip\noindent
{\bf 3.6.~Lemma.}
{\it With the notation of {\rm 3.1}, assume
$ \Gr_{n-1+k}^{W}H^{n-1}(F_{x},\bC)_{\lambda} \ne 0 $ for a positive
integer
$ k $,
where
$ W $ is the weight filtration.
Then
$ N^{k} \ne 0 $ on
$ \psi_{f,\lambda}\bC_{X} $,
where
$ N $ is the logarithm of the unipotent part of the monodromy.
}

\medskip\noindent
{\it Proof.}
We have the weight filtration
$ W $ on the perverse sheaf
$ \psi_{f,\lambda}\bC_{X}[n-1] $ (see [2]) such that we
have isomorphisms for
$ j > 0 $
$$
N^{j} : \Gr_{n-1+j}^{W}\psi_{f,\lambda}\bC_{X}[n-1]
\to \Gr_{n-1-j}^{W}\psi_{f,\lambda}\bC_{X}[n-1].
$$
This gives the weight filtration of a mixed Hodge module, see
[33], [34].
Furthermore, the mixed Hodge structure on
$ H^{n-1}(F_{x},\bC) $ is given by applying the pull-back functor
$ H^{j}{i}_{x}^{*} $ to
$ \psi_{f,\lambda}\bC_{X}[n-1] $ as in 3.5.
The functor
$ H^{0}{i}_{x}^{*} $ preserves the condition that
$ \Gr_{i}^{W} = 0 $ for
$ i > r $ where
$ r $ is any fixed integer, see [34].
So the hypothesis implies that
$ \Gr_{n-1+m}^{W}\psi_{f,\lambda}\bC_{X}[n-1] \ne 0 $ for some
$ m \ge k $,
and the assertion follows.

\bigskip\bigskip
\centerline{{\bf 4. Brieskorn modules and Gauss-Manin systems}}

\bigskip\noindent
{\bf 4.1.}
Let
$ f $ be a non-constant holomorphic function on a complex manifold
$ X $ of dimension
$ n \ge 2 $, and
$ x \in D := f^{-1}(0) $.
With the notation of 3.1, the Gauss-Manin system is defined by
$$
\cG_{f} = H^{0}\DR_{X}\cB_{f,x}\,(= \omega_{X,x}
\otimes_{\cD_{X,x}}\cB_{f,x}),
\leqno(4.1.1)
$$
where
$ \omega_{X} $ is the sheaf of the differential
forms of degree
$ n $.
Here we consider only the cohomology of degree
$ 0 $ because we assume essentially isolated singularity conditions
when we consider Gauss-Manin systems in this paper.
This is a regular holonomic
$ \bC\{t\}\langle \partial_{t} \rangle $-module where
$ \bC\{t\}\langle \partial_{t} \rangle = \cD_{S,0} $ with
$ S $ an open disc.
It is known that
$ \cG_{f} $ is a finite free
$ \bC\{\!\{\partial_{t}^{-1}\}\!\}[\partial_{t}] $-module of rank
$ \mu_{n-1} $ where
$ \mu_{j} $ is the rank of the
$ j $-th cohomology of the Milnor fiber around
$ x $, see e.g. [1].

The Brieskorn module is defined by
$$
\cH''_{f} = \Omega_{X,x}^{n}/df\wedge d\Omega_{X,x}^{n-2}.
$$
It is a
$ \bC\{t\}\langle\partial_{t}^{-1}\rangle $-module, where
$ \partial_{t}^{-1}\omega $ is defined by
$ df\wedge\eta $ with
$ d\eta = \omega $.
There is a canonical morphism
$$
\cH''_{f} \to \cG_{f},
\leqno(4.1.2)
$$
compatible with the action of
$ \bC\{t\}\langle\partial_{t}^{-1}\rangle $
so that
$ \cG_{f} $ is identified with the localization of
$ \cH''_{f} $ by
$ \partial_{t}^{-1} $.

Let
$ V $ be the filtration of Kashiwara and Malgrange on
$ \cG_{f} $ indexed by
$ \bQ $ so that
$$
\partial_{t}t - \alpha \,\,\,\text{is nilpotent on}\,\,\,
\Gr_{V}^{\alpha}\cG_{f}.
$$
It is known that
$ V^{\alpha}\cG_{f} $ for
$ \alpha > 0 $ is naturally identified with the Deligne extension
of the restriction to a punctured disk of a coherent extension of
$ \cG_{f} $ such that the eigenvalues of the
residue of the connection are contained in
$ [\alpha-1,\alpha) $.
In particular, we have isomorphisms for
$ \alpha \in (0,1] $ and
$ \lambda = \exp(-2\pi i\alpha) $
$$
H^{n-1}(F_{x},\bC)_{\lambda}\simeq\Gr_{V}^{\alpha}\cG_{f},
\leqno(4.1.3)
$$
where the left-hand side is the
$ \lambda $-eigenspace of the Milnor cohomology for the Milnor
monodromy.
We define the Brieskorn lattices of
$ \cG_{f} $ by
$$
\cG_{f}^{(i)} = \partial_{t}^{-i}\cG_{f}^{(0)}\,(i \in \bZ)
\quad\text{with}\quad
\cG_{f}^{(0)} = \hbox{\rm Im}(\cH''_{f} \to \cG_{f}),
$$
where the last morphism
$ \cH''_{f} \to \cG_{f} $ is as in (4.1.2).
Note that
$ \cG_{f}^{(i)} \subset \cG_{f}^{(i-1)} $ because
$ \cG_{f}^{(0)} $ is stable by the action of
$ \partial_{t}^{-1} $.
Let
$ \tcG_{f}^{(i)} $ be the saturations of
$ \cG_{f}^{(i)} $, i.e.
$$
\tcG_{f}^{(i)} =
\msum_{k\ge 0}(\partial_{t}t)^{k}\cG_{f}^{(i)} =
\msum_{k\ge 0}(t\partial_{t})^{k}\cG_{f}^{(i)}.
$$
They have the induced filtration
$ V $.
By [17], [23], we have
$$
\tcG_{f}^{(0)} \subset V^{>0}\cG_{f}.
\leqno(4.1.4)
$$
This implies
$ \partial_{t}t\tcG_{f}^{(i)} = \tcG_{f}^{(i)} $ for
$ i \ge 0 $, and
$$
\tcG_{f}^{(i)} = \partial_{t}^{-i} \tcG_{f}^{(0)} = t^{i}
\tcG_{f}^{(0)}\quad\text{for}\,\, i \in \bN.
$$
For
$ \alpha > 0 $ and
$ i \in \bN $, we have
$$
t^{i} : \Gr_{V}^{\alpha}\tcG_{f}^{(-i)}\buildrel\sim\over\to
\Gr_{V}^{\alpha+i}\tcG_{f}^{(0)}\subset\Gr_{V}^{\alpha+i}\cG_{f},
\leqno(4.1.5)
$$
because the action of
$ t^{i}\partial_{t}^{i} = t\partial_{t}(t\partial_{t}-1)\cdots
(t\partial_{t}-i+1) $ on
$ \Gr_{V}^{\alpha+i}\tcG_{f}^{(0)} $ for
$ \alpha > 0 $ is injective and hence surjective.

Using the isomorphism (4.1.3) for
$ \alpha \in (0,1] $ and
$ \lambda = \exp(-2\pi i\alpha) $, we define decreasing filtrations
$ P $ and
$ \tP $ on the Milnor cohomology by
$$
\aligned
P^{n-1-i} H^{n-1}(F_{x},\bC)_{\lambda}
&\simeq\Gr_{V}^{\alpha}\cG_{f}^{(-i)},
\\
\tP^{n-1-i} H^{n-1}(F_{x},\bC)_{\lambda}
&\simeq\Gr_{V}^{\alpha}\tcG_{f}^{(-i)}.
\endaligned
\leqno(4.1.6)
$$
Note that
$ P^{n-1-i} = \tP^{n-1-i} = 0 $ for
$ i < 0 $ by (4.1.4), and
$ \tP^{n-1-i} $ is stable by the Milnor monodromy because
$ \tcG_{f}^{(-i)} $ is stable by the action of
$ t\partial_{t} $.

If
$ \xi f =f $ for a vector field
$ \xi $, then
$ \cG_{f}^{(-i)} = \tcG_{f}^{(-i)} $ and
$ P = \tP $.
In the isolated singularity case,
$ P $ coincides with the Hodge filtration, see [36], [42].
Note that
$ \cG_{f}^{(1)} $ is the image of
$ df\wedge\Omega_{X,x}^{n-1} $ and
$ \cG_{f}^{(0)}/\cG_{f}^{(1)} $ is a quotient of
$$
\Omega_{f} := \omega_{X,x}/df\wedge {\Omega}_{X,x}^{n-1},
$$
because
$ d(\eta\otimes 1) = d\eta\otimes 1 - (df\wedge\eta)\otimes
\partial_{t} $ where
$ d $ is the differential of
$ \DR_{X}\cB_{f} $, see (3.1.1).
If
$ D $ has an isolated singularity, then it is well known that
(4.1.2) is injective, and
$ \cG_{f}^{(0)}/\cG_{f}^{(1)} = \Omega_{f} $, see [29], [36].

\medskip\noindent
{\bf 4.2.~Proposition.}
{\it The filtration
$ V $ on the Gauss-Manin system
$ \cG_{f} $ coincides with the filtration induced by the filtration
$ V $ of Kashiwara and Malgrange on
$ \cB_{f} $ via the isomorphism {\rm (4.1.1)}
using any trivialization of
$ \omega_{X,x} $.
We have the canonical isomorphisms for
$ \alpha \in \bQ $
}
$$
V^{\alpha}\cG_{f} = H^{0}\DR_{X}(V^{\alpha}\cB_{f,x}),\quad
\Gr_{V}^{\alpha}\cG_{f} = H^{0}\DR_{X}(\Gr_{V}^{\alpha}\cB_{f,x}).
\leqno(4.2.1)
$$

\medskip\noindent
{\it Proof.}
For
$ -\infty < \alpha < \beta < \gamma < +\infty $, we have a long
exact sequence
$$
\buildrel{\partial}\over\to
H^{i}\DR_{X}(V^{\beta}/V^{\gamma})\cB_{f,x} \to
H^{i}\DR_{X}(V^{\alpha}/V^{\gamma})\cB_{f,x} \to
H^{i}\DR_{X}(V^{\alpha}/V^{\beta})\cB_{f,x}
\buildrel{\partial}\over\to,
$$
compatible with the action of
$ \partial_{t}t $ so that the connecting morphisms
$ \partial $ vanish, where
$ (V^{\alpha}/V^{\beta})\cB_{f,x} =
V^{\alpha}\cB_{f,x}/V^{\beta}\cB_{f,x} $.
We have the finiteness of
$ H^{i}\DR_{X}(V^{\alpha}\cB_{f,x}) $ over
$ \bC\{t\} $ by [33], 3.4.8, and the connecting morphisms
$ \partial $ vanish also for
$ \gamma = +\infty $ (where
$ V^{+\infty} = 0 $) using the completion, see loc.~cit.
So the first isomorphism follows, and the second isomorphism then
follows using the vanishing of
$ \partial $.

\medskip\noindent
{\bf 4.3.~Proposition.}
{\it With the notation of {\rm 3.1} we have for any
$ \alpha \in \bQ $
}
$$
\aligned
\hbox{Im}(H^{0}\DR_{X}(V^{\alpha}M_{x})\to H^{0}\DR_{X}(\cB_{f,x}))
&=V^{\alpha}\tcG_{f}^{(0)},
\\
\hbox{Im}(H^{0}\DR_{X}(\Gr_{V}^{\alpha}M_{x})\to
H^{0}\DR_{X}(\Gr_{V}^{\alpha}\cB_{f,x}))
&=\Gr_{V}^{\alpha}\tcG_{f}^{(0)}.
\endaligned
$$

\medskip\noindent
{\it Proof.}
The canonical morphism
$ H^{0}\DR_{X}(V^{\alpha}M_{x}) \to
H^{0}\DR_{X}(\Gr_{V}^{\alpha}M_{x}) $ is surjective, because
$ H^{1}\DR_{X}(V^{>\alpha}M_{x}) = 0 $.
So the second isomorphism is reduced to the first using Proposition~4.2,
and it is enough to show the first isomorphism.
The right-hand side is the intersection of the images of
$ H^{0}\DR_{X}(M_{x}) $ and
$ H^{0}\DR_{X}(V^{\alpha}\cB_{f,x}) $, because
$ \tcG_{f}^{(0)} $ is the image of
$ H^{0}\DR_{X}(M_{x}) $.
By the commutative diagram
$$
\CD
H^{0}\DR_{X}(V^{\alpha}M) @>>> H^{0}\DR_{X}(V^{\alpha}\cB_{f})
@>>> H^{0}\DR_{X}V^{\alpha}(\cB_{f}/M) @>>> 0
\\
@VVV @VVV @VVV
\\
H^{0}\DR_{X}(M) @>>> H^{0}\DR_{X}(\cB_{f}) @>>>
H^{0}\DR_{X}(\cB_{f}/M) @>>> {\,0,}
\endCD
$$
the assertion is reduced to the injectivity of the last vertical
morphism, but this is easily proved by using the action of
$ \partial_{t}t $ together with a long exact sequence as in the proof
of Proposition~4.2.
So the assertion follows.

\medskip\noindent
{\bf 4.4.~Proposition.}
{\it With the notation of {\rm 4.1} the Hodge filtration
$ F $ on the Milnor cohomology is contained in
$ \tP $.
}

\medskip\noindent
{\it Proof.}
The Hodge filtration
$ F $ on the Milnor cohomology is defined by using the
construction in 3.5.
For any regular holonomic
$ \cD_{X^{j-1}} $-module
$ N $, we have canonical morphisms of complexes
$$
C(\partial_{j}:\Gr_{V_{j}}^{1}N\to \Gr_{V_{j}}^{0}N)
\leftarrow
C(\partial_{j}:V_{j}^{1}N\to V_{j}^{0}N)\to
C(\partial_{j}:N\to N),
$$
which are quasi-isomorphisms at least after taking
the de Rham functor on
$ X^{j} $.
Iterating this, we get a canonical isomorphism in the derived
category
$$
i_{x}^{*}\Gr_{V}^{\alpha}\cB_{f}=
\DR_{X}(\Gr_{V}^{\alpha}\cB_{f,x}),
$$
where the left-hand side is defined as in 3.5.
If
$ N $ underlies a mixed Hodge module so that it has the
Hodge filtration
$ F $, then the Hodge filtration
$ F $ on
$ \Gr_{V_{j}}^{\alpha}N $ is induced by
$ F $ on
$ N $, and the canonical surjection
$ V_{j}^{0}N \to \Gr_{V_{j}}^{0}N $ is strictly compatible with
$ F $.
This implies that the Hodge filtration
$ F $ on 
$ H^{0}i_{x}^{*}\Gr_{V}^{\alpha}\cB_{f} $ is contained in the
filtration on
$ H^{0}\DR_{X}(\Gr_{V}^{\alpha}\cB_{f,x}) =
\Gr_{V}^{\alpha}\cG_{f} $ induced by
$ F $ on
$ \Gr_{V}^{\alpha}\cB_{f,x} $ (up to an appropriate shift)
via the above isomorphism.
Moreover, the latter filtration on
$ \Gr_{V}^{\alpha}\cG_{f} $ is contained in
$ \tP $ by (4.1.5) and Proposition~4.3.
Indeed,
$$
F_{p}\cB_{f,x} = \moplus_{0\le i\le p}\cO_{X}\otimes
\partial_{t}^{i},
$$
and the image of
$ F_{p}V^{\alpha}\cB_{f,x} $ in
$ H^{0}\DR_{X}(\cB_{f,x}) = \cG_{f} $ is contained in
$ \tcG^{(-p)} \cap V^{\alpha}\cG_{f} $ by
$ \tcG_{f}^{(-i)}\subset\tcG_{f}^{(-i-1)} $.
So the assertion follows.

\medskip\noindent
{\bf 4.5.~Proof of Theorem~2.}
Since the de Rham complex
$ \DR_{X} $ is the Koszul complex for
$ \partial_{1},\dots,\partial_{n} $
(trivializing
$ \omega_{X} $ by
$ dx_{1}\wedge\cdots\wedge dx_{n} $) and
$ M_{x} $ is generated over
$ \cD_{X}[s] $ by
$ 1\otimes 1 $ (or
$ f^{s} $), we see that the image of
$ H^{0}\DR_{X}M_{x} $ in
$ \cG_{f} $ coincides with
$ \tcG_{f}^{(0)} $.
By (4.1.3) the filtration
$ \tP $ in (4.1.6) is identified with a filtration on
$ \Gr_{V}^{\alpha}\cG_{f} $ for
$ \alpha \in (0,1] $, which is also denoted by
$ \tP $.
Then, by (4.1.5), we have for
$ i \in \bN $
$$
t^{i} : \tP^{n-1-i}\Gr_{V}^{\alpha}\cG_{f}\buildrel\sim\over\to
\Gr_{V}^{\alpha+i}\tcG_{f}^{(0)}\subset\Gr_{V}^{\alpha+i}\cG_{f}.
$$

Let
$ \tP $ denote also the filtration on
$ \Gr_{V}^{\alpha}\cB_{f} $ for
$ \alpha \in (0,1] $ such that for
$ i \in \bN $
$$
t^{i} : \tP^{n-1-i}\Gr_{V}^{\alpha}\cB_{f}\buildrel\sim\over\to
\Gr_{V}^{\alpha+i}M\subset\Gr_{V}^{\alpha+i}\cB_{f}.
$$
Then
$$
t^{i} : \Gr_{\tP}^{n-1-i}\Gr_{V}^{\alpha}\cB_{f}
\buildrel\sim\over\to\Gr_{V}^{\alpha+i}(M/tM).
$$
By Proposition~4.3 the filtration
$ \tP $ on
$ \Gr_{V}^{\alpha}\cB_{f} $ induces
$ \tP $ on
$ \Gr_{V}^{\alpha}\cG_{f} $ taking the de Rham functor
$ \DR_{X} $.
So the assertion
$ (a) $ follows from the spectral sequence
$$
E_{1}^{p,q}=H^{p+q}\DR_{X}\Gr_{\tP}^{p}\Gr_{V}^{\alpha}\cB_{f,x}
\Rightarrow H^{p+q}\DR_{X}\Gr_{V}^{\alpha}\cB_{f,x},
\leqno(4.5.1)
$$
because
$ E_{1}^{p,q} \ne 0 $ if
$ E_{\infty}^{p,q} \ne 0 $.

If
$ \alpha + i $ is not a root of
$ b_{f,y}(-s) $ for any
$ y \ne x $ and
$ i \ge i_{0} $ where
$ i_{0} $ is a nonnegative integer, then
$ E_{1}^{p,q} = 0 $ for
$ p + q > 0 $ or
$ p + q < 0 $ and
$ p \le p_{0} := n - 1 - i_{0} $.
Indeed, a
$ \cD_{X} $-module supported on a point is a direct sum of
$ \bC[\partial] $ in the notation of 3.3, and
$ \DR_{X}(\bC[\partial]) = \bC $.
So we have
$ E_{1}^{p,q} = E_{\infty}^{p,q} $ for
$ p \le p_{0} $, and the assertion
$ (b) $ follows.
If the assumption of
$ (c) $ is satisfied, then
$ E_{1}^{p,q} = 0 $ for
$ p + q \ne 0 $ and
$ E_{1}^{p,q} = E_{\infty}^{p,q} $ for any
$ p, q $ by a similar argument.
So the assertion follows.

\medskip\noindent
{\bf 4.6.~Proposition.}
{\it With the notation of {\rm 3.1} and {\rm 4.1}, assume
$ f $ is a weighted homogeneous polynomial of strictly positive
weights
$ (w_{1},\dots,w_{n}) $.
Let
$ \xi = \sum_{i}w_{i}x_{i}\partial_{i} $
so that
$ \xi f = f $ and hence
$ P = \tP $.
Let
$ \Omega_{f}^{\beta} $ be the
$ \beta $-eigenspace for the Lie derivation by
$ \xi $.
Let
$ \alpha $ be a rational number such that
$ \exp(-2\pi i\alpha) $ is not an eigenvalue of the Milnor
monodromy at
$ y \ne x $.
Then
$ \alpha $ is a root of
$ b_{f}(-s) $ if and only if the image of
$ \Omega_{f}^{\alpha} $ in
$ \cG_{f}^{(0)}/\cG_{f}^{(1)} $ does not vanish.
Its multiplicity is
$ 1 $ if
$ \alpha $ is a root.
}

\medskip\noindent
{\it Proof.}
Since
$ \xi f = f $, we have
$ t\cG_{f}^{(0)} \subset \partial_{t}^{-1}\cG_{f}^{(0)} $ and hence
$ \tcG_{f}^{(0)} = \cG_{f}^{(0)} $ by the definition of the action of
$ \partial_{t}^{-1} $.
Let
$ {\omega}_{X,x}^{\beta} $ be the
$ \beta $-eigenspace for the action of
$ \xi $.
Then the action of
$ \partial_{t}t $ on the image of
$ {\omega}_{X,x}^{\beta}\otimes 1\otimes 1 $ in
$ \cG_{f} $ is given by the multiplication by
$ \beta $.
Indeed, if we take
$ gdx_{1}\wedge\cdots\wedge dx_{n} \in {\omega}_{X,x}^{\beta} $
so that
$ \xi g = (\beta-\sum_{i}w_{i})g $, then we have by (3.1.1)
$$
(\msum_{i}w_{i}\partial_{i}x_{i})(g\otimes 1)=
\beta g\otimes 1 - \partial_{t}t(g\otimes 1),
$$
where the left-hand side vanishes in the cohomology of the
de Rham complex which is identified with the Koszul complex for
$ \partial_{1},\dots,\partial_{n} $ as above, see also [32].

This implies that
$ \Gr_{V}^{\beta}\cG_{f}^{(0)} $ is generated by the image of
$ \omega_{X,x}^{\beta} $ so that the action of
$ \partial_{t}t $ on
$ \Gr_{V}^{\beta}\cG_{f}^{(0)} $ is semisimple
(using the algebraic Gauss-Manin system if necessary).
Then the assertions follow from Theorem~2.

\medskip
In the isolated singularity case, we have the following:

\medskip\noindent
{\bf 4.7.~Proposition.}
{\it With the notation of {\rm 4.1}, assume
$ D $ has an isolated singularity at
$ x $.
Let
$ V $ denote the filtration on
$ \omega_{X,x} = \Omega_{X,x}^{n} $ induced by the filtration
$ V $ on
$ \cO_{X,x} $ using any trivialization of
$ \omega_{X,x} $.
Then the natural projection
$ \omega_{X,x}\to \cH''_{f} $ is strictly compatible
with the filtration
$ V^{\alpha} $ for
$ \alpha \le 1 $.
}

\medskip\noindent
{\it Proof.}
Let
$ V' $ denote the filtration on
$ \omega_{X,x} $ induced by the filtration
$ V $ on
$ \cH''_{f} $ using the projection
$ \omega_{X,x}\to \cH''_{f} $.
Since the filtration
$ V $ on
$ \cB_{f} $ induces the filtration
$ V $ on the Gauss-Manin system by Proposition~4.2, we have
$ V^{\alpha}\omega_{X,x} \subset V'{}^{\alpha}\omega_{X,x} $.
Then we get the equality for
$ \alpha \le 1 $ by calculating the dimension of their graded
pieces for
$ \alpha < 1 $, because they both give the coefficient
$ m_{\alpha} $ of the spectrum for
$ \alpha < 1 $.
So the assertion follows.

\medskip\noindent
{\bf 4.8.~Hodge and pole order filtrations.}
Assume
$ X $ is affine space
$ \bC^{n} $, and
$ D $ is the affine cone of a divisor
$ Z $ of degree
$ d $ on
$ Y := \bP^{n-1} $.
Then there is a cyclic covering
$ \pi : \tY \to Y $ of degree
$ d $ which is ramified along
$ Z $.
Put
$ U = Y \setminus Z $,
$ \tU = \pi^{-1}(U) $.
Then
$ \tU $ is identified with the Milnor fiber
$ F_{0} $ of a function
$ f $ defining the affine cone
$ D $ of
$ Z $, and the geometric Milnor monodromy corresponds to a generator
of the covering transformation group of
$ \tU \to U $, see [11], 1.8.

For
$ k = 1, \dots, d $, let
$ L^{(k/d)} $ be the direct factor of
$ \pi_{*}\bC_{\tU} $ on which the action of the Milnor monodromy is the
multiplication by
$ \exp(-2\pi ik/d) $ so that
$ H^{j}(U,L^{(k/d)}) = H^{j}(F_{0},\bC)_{\lambda} $ where
$ \lambda = \exp(-2\pi ik/d) $.
Note that
$ L^{(k/d)} $ is a local system of rank
$ 1 $ on
$ U $, and its monodromy around a smooth point
$ x $ of
$ Z_{\red} $ is the multiplication by
$ \exp(2\pi ikm_{x}/d) $ where
$ m_{x} $ is the multiplicity of
$ Z $ at
$ x $.
(This can be shown by blowing up along the origin of
$ \bC^{n} $ and considering the nearby cycles for the pull-back of
$ f $, see also [10].)

Let
$ \cL^{(k/d)} $ be the meromorphic extension of
$ L^{(k/d)}\otimes_{\bC}\cO_{U} $.
This is a regular holonomic
$ \cD_{Y} $-module on which the action of a function
$ h $ defining
$ Z $ is bijective.
We see that
$ \cL^{(k/d)} $ is locally isomorphic to a free
$ \cO_{Y}(*Z) $-module generated by a multivalued function
$ h_{j}^{-k/d} $ where
$ h_{j} = x_{j}^{-d}f $ on
$ \{x_{j} \ne 0 \} \subset \bP^{n-1} $.
Note that the
$ \cO_{Y} $-submodule generated locally by
$ h_{j}^{-k/d}\,(= x_{j}^{k}f^{-k/d}) $ is isomorphic to
$ \cO_{Y}(k) $, because the relation
$ g_{j}h_{j}^{-k/d} = g_{i}h_{i}^{-k/d} $ means that
$ g_{j}x_{j}^{k} = g_{i}x_{i}^{k} $, i.e.
$ \{g_{j}\} $ defines a section of
$ \cO_{Y}(k) $.

The pole order filtration
$ P_{i}\cL^{(k/d)} $ is defined to be the locally free
$ \cO_{Y} $-submodule of
$ \cL^{(k/d)} $ generated by
$ h_{j}^{-i-(k/d)} $ on
$ \{x_{j} \ne 0 \} $ for
$ i \in \bN $, and
$ P_{i}\cL^{(k/d)} = 0 $ for
$ i < 0 $.
Then
$ P_{i}\cL^{(k/d)} $ is isomorphic to
$ \cO_{Y}(id+k) $ by the above argument.
On the other hand, there is the Hodge filtration
$ F $ on
$ \cL^{(k/d)} $ such that
$ F_{i}\cL^{(k/d)} = P_{i}\cL^{(k/d)} $ outside
$ \Sing Z_{\red} $ for any
$ i $ by the theory of mixed Hodge modules.
Then we have
$ F_{i}\cL^{(k/d)} \subset P_{i}\cL^{(k/d)} $ on
$ Y $ because
$ P_{i}\cL^{(k/d)} $ is locally free and
$ \Sing Z_{\red} $ has codimension
$ \ge 2 $ in
$ Y $.

The Hodge and pole order filtrations are closely related respectively
to the spectrum and the
$ b $-function of
$ f $.
Indeed, the Hodge filtration
$ F $ on
$ \cL^{(k/d)} $ induces the Hodge filtration on the Milnor cohomology
by taking the de Rham cohomology.
Similarly the pole order filtration
$ P $ on the Milnor cohomology is defined by using the de Rham
cohomology.
Here the filtration is shifted by the degree of the differential forms,
and the associated decreasing filtration is used.
Then we have Theorem~2 together with the following

\medskip\noindent
{\bf 4.9.~Proposition.}
{\it With the above notation and assumption, the above pole order
filtration
$ P $ coincides with the filtration
$ \tP = P $ in {\rm (4.1.6)}.
Moreover, for
$ \alpha = k/d \in (0,1) $ and
$ \lambda = \exp(-2\pi i\alpha) $, the above
$ P^{n-1-i} $ is identified with the image of
$$
\Gr_{V}^{\alpha+i}\cG_{f}^{(0)} \subset \Gr_{V}^{\alpha+i}\cG_{f}
\simeq\Gr_{V}^{\alpha}\cG_{f} \simeq H^{n-1}(F_{0},\bC)_{\lambda},
$$
where the middle isomorphism can be induced by both
$ \partial_{t}^{k} $ and
$ t^{-k} $, and the last morphism is induced by {\rm (4.1.3)}.
}

\medskip\noindent
{\it Proof.}
This follows from the arguments in [11], using the local
generator
$ h_{j}^{-k/d} $ in 4.8 to define an isomorphism generalizing
Lemma~1.2 in loc.~cit.
Note that
$ \cS^{i} \,(0 < i < d) $ in loc.~cit. is identified with
$ \cL^{((d-i)/d)}\otimes \cO_{Y}(-Z) $ in this paper, and
$ \Omega^{j}[f^{-1}]_{k}^{(\xi)} $ (the degree
$ k $ part of the image of the interior product
$ \iota_{\xi} $) is identified with the vector space of
meromorphic sections of
$ \pi_{*}\pi^{*}\Omega_{Y}^{j} $ over
$ Y\setminus Z $ on which the Lie derivation
$ L_{\xi} $ acts as the multiplication by
$ k $.
Here
$ \pi : \bC^{n}\setminus \{0\} \to Y $ denotes the canonical
projection.
Then we get the desired isomorphism by using the restriction to
$ \{x_{i}=1\}\subset\bC^{n} $ for any
$ i $.
The above identification of the filtrations is compatible with
(4.1.6) because
$ t^{i}\Gr_{V}^{\alpha}\partial_{t}^{i}\cG_{f}^{(0)}
= t^{i}\partial_{t}^{i}\Gr_{V}^{\alpha+i}\cG_{f}^{(0)}
= \Gr_{V}^{\alpha+i}\cG_{f}^{(0)} $.

\medskip\noindent
{\bf 4.10.~Remark.}
If
$ Z $ is smooth, the two filtrations
$ F_{i} $ and
$ P_{i} $ on
$ \cL^{(k/d)} $ coincide for any
$ i $, and this explains the coincidence of the spectrum and the
roots of
$ b_{f}(-s) $ (forgetting the multiplicity) in this case.
However, if
$ Z $ is a reduced divisor with normal crossings, these two
filtrations coincide only for
$ i = 0 $, and not for
$ i > 0 $ because the Hodge filtration is defined by using the
{\it sum} of the pole orders along the irreducible components,
see [7].
This explains the fact that the spectrum and the roots of
$ b_{f}(-s) $ coincide (forgetting the multiplicity)
only if they are restricted to the interval
$ (0,1] $ in this case.

\medskip\noindent
{\bf 4.11.~Remark.}
Let
$ D $ be as in 4.8 so that
$ f $ is a homogeneous polynomial of degree
$ d $.
Then the Brieskorn lattice
$ \cG_f^{(0)} $ has a monomial basis
$ (\omega_j) $ over
$ \bC\{t\} $ such that each
$ \omega_j $ is represented by
$$
x^{\nu}dx_1\wedge\cdots\wedge dx_n \quad\text{with}\,\,\,
\nu=(\nu_1,\dots,\nu_n)\in\bN^n.
$$
Restricting to a Milnor fiber, this gives a basis
$ ([\omega_j]) $ of the
Milnor cohomology.
For the last assertion, it is enough to assume that
$ (\omega_j) $ gives a basis of
$ \cG_f^{(0)}[t^{-1}] $ over
$ \bC\{t\}[t^{-1}] $, and the minimality condition in
[43], Thm.~4.8 corresponds to that it gives a basis of
$ \cG_f^{(0)}$ over
$ \bC\{t\} $.
Using this basis, the pole order filtration
$ P $ on the Milnor cohomology is expressed as
$$
P^{n-1-i} = \msum_{\deg(\omega_j)/d\le i+1}\bC[\omega_j],
$$
where
$ \deg(x^{\nu}dx_1\wedge\cdots\wedge dx_n)=|\nu|+n $.
This follows from the well-known formula
$$
\partial_tt(x^{\nu}dx_1\wedge\cdots\wedge dx_n)=
(\deg(x^{\nu}dx_1\wedge\cdots\wedge dx_n)/d)
x^{\nu}dx_1\wedge\cdots\wedge dx_n.
$$
(See the proof of Prop.~4.6.)
So Theorem~2 (a) may be viewed as a generalization of
[43], Thm.~4.8.

\bigskip\bigskip
\centerline{\bf 5. Case of hyperplane arrangements}

\bigskip\noindent
{\bf 5.1.~Cohomology of twisted de Rham complexes.}
With the notation and the assumptions of 4.8, assume further
that
$ D $ is the affine cone of a projective hyperplane arrangement
$ Z $ in
$ Y=\bP^{n-1} $, i.e.
$ D $ is a central hyperplane arrangement.
Then, by [3], [14], [35], the cohomology of the local systems on
$ U = Y \setminus Z $ in 4.8 can be calculated as follows:

Let
$ Z_{i} \,(1 \le i \le d) $ be the irreducible components of
$ Z $ where
$ d = \deg Z $, and
$ x_{1},\dots,x_{n} $ be coordinates of
$ \bC^{n} $ such that
$ Z_{d} = \{x_{n}=0\} $.
Then the complement
$ Y' $ of
$ Z_{d} $ in
$ Y $ is identified with
$ \bC^{n-1} $.
Let
$ g_{i} $ be a polynomial of degree
$ 1 $ on
$ Y' $ defining
$ Z_{i} \cap Y' $.
Put
$$
\omega_{i}= dg_{i}/g_{i} \,(1 \le i \le d-1),\quad
h = g_{1}\cdots g_{d-1}.
$$
For
$ \alpha = (\alpha_{1},\dots,\alpha_{d-1}) \in \bC^{d-1} $, let
$$
h^{\alpha} = g_{1}^{\alpha_{1}}\cdots g_{d-1}^{\alpha_{d-1}},
$$
and
$ \cO_{Y'}h^{\alpha} $ be a free
$ \cO_{Y'} $-module of rank
$ 1 $ on
$ Y' $ with formal generator
$ h^{\alpha} $.
There is a regular singular integrable connection
$ \nabla $ such that for
$ u\in\cO_{Y'} $
$$
\nabla(uh^{\alpha}) = (du)h^{\alpha} + u\omega h^{\alpha}
\quad\text{with}\quad
\omega=\msum_{1\le i\le d-1}\alpha_{i}\omega_{i}.
$$
Let
$ \cA^{p}_{h,\alpha} $ be the
$ \bC $-vector subspace of
$ \Gamma(U,\Omega_{U}^{p}h^{\alpha}) $ generated by
$ \omega_{i_{1}}\wedge\cdots\wedge\omega_{i_{p}}h^{\alpha} $ for any
$ i_{1}<\cdots<i_{p} $.
Then
$ \cA^{\ssbull}_{h,\alpha} $ with differential
$ \omega\wedge $ is a subcomplex of
$ \Gamma(U,\Omega_{U}^{\ssbull}h^{\alpha}) $.
Put
$ \alpha_{d} = -\msum_{1\le i\le d-1}\alpha_{i} $.
By [3], [14], [35], we have the canonical quasi-isomorphism
$$
\cA^{\ssbull}_{h,\alpha} \buildrel\sim\over\longrightarrow
\Gamma(U,\Omega_{U}^{\ssbull}h^{\alpha}),
\leqno(5.1.1)
$$
if the following condition holds for any dense edge
$ L $ of
$ Z $:
$$
\alpha_{L} := \msum_{Z_{i}\supset L}\alpha_{i} \notin
\bN\setminus \{0\}.
\leqno(5.1.2)
$$
Here an edge is an intersection of
$ Z_{i} $ over a subset of
$ \{1,\dots,d\} $, and
an edge is called dense if and only if
the hyperplanes containing the edge are identified with an
indecomposable central arrangement
(where an arrangement in
$ \bC^{n} $ is called decomposable if and only if there is a
decomposition
$ \bC^{n} = \bC^{n'}\times\bC^{n''} $ such that the arrangement is
the union of the pull-backs of arrangements on
$ \bC^{n'} $ and
$ \bC^{n''} $, see [35] for details.)
In the case of a constant local system, this is due to [3].
In a general case it is shown in [14],
and is improved in [35].
Note that if
$ Z $ is a divisor with normal crossings (i.e. if
$ Z $ is generic), then condition (5.1.2) is equivalent to
$ \alpha_{i} \notin \bN\setminus \{0\} $ for any
$ i \in [1,d] $ (because the dense edges
consist of the
$ Z_{i} $ in this case), and [14] is sufficient in this case.

As a corollary, we get

\medskip\noindent
{\bf 5.2.~Proposition.}
{\it Let
$ D $ be a central hyperplane arrangement in
$ \bC^{n} $ defined by a reduced polynomial
$ f $ of degree
$ d $.
Let
$ Z $ be the projective arrangement in
$\bP^{n-1} $ corresponding to
$ D $.
For an edge
$ L $ of
$ Z $, let
$ m_{L} $ be the number of hyperplanes
$ Z_{i} \supset L $.
Assume all the roots of
$ b_{h,y}(-s) $ are strictly less than
$ 2 $ for any
$ y \in Z $ where
$ h $ is a reduced local equation of
$ Z $ at
$ y $.
Assume, moreover, there is a hyperplane, denoted by
$ Z_{d} $, such that
$ \hbox{\rm GCD}(m_{L},d) = 1 $ for any dense edge
$ L $ of
$ Z $ contained in
$ Z_{d} $.
Then the roots of
$ b_{f,0}(-s) $ are strictly less than
$ 2 $.
}

\medskip\noindent
{\it Proof.}
We apply the above argument to the case
$ \alpha_{i} = -k/d \,\, (0 \le i < d) $ for each
$ k \in [0,d-1] $.
Then for
$ k \in [1,d-1] $, we see that
$ h^{\alpha} = h^{-k/d} $ is a section of
$ P_{0}\cL^{(k/d)} $ which has a zero of order
$ k $ along the divisor at infinity
$ Z_{d} $, and condition (5.1.1) is satisfied for any dense edge
$ L $ of
$ Z $.
(Indeed,
$ \alpha_{L} \notin \bZ $ if
$ L \subset Z_{d} $, and
$ \alpha_{L} \le 0 $ otherwise.)
Moreover, the meromorphic extension of
$ \omega_{i_{1}}\wedge\cdots\wedge\omega_{i_{p}} $ to
$ Y $ has at most a pole of order
$ 1 $ along each
$ Z_{i} \,\,(1 \le i \le d) $.
Thus
$ \cA_{h,\alpha}^{n-1} $ is contained in
$ \Omega_{Y}^{n-1}\otimes_{\cO}P_{i}\cL^{(k/d)} $ with
$ i = 1 $ if
$ 0 < k < d $ and
$ i = 0 $ if
$ k = 0 $ by the definition of
$ P_{i} $.
So the assertion follows from Theorem~2 together with (5.1.1)
using [3], [14], [35].

\medskip
The first assertion of Proposition below is due to [43].

\medskip\noindent
{\bf 5.3.~Proposition.}
{\it With the above notation,
$ -1 $ is the only integral root of
$ b_{f}(s) $ {\rm (}see {\rm [43])}, and its multiplicity is
$ n $, assuming the arrangement is not the pull-back of an
arrangement in a strictly lower dimensional vector space.
}

\medskip\noindent
{\it Proof.}
The assertion is well known in the normal crossing case.
In particular, it holds on the smooth part of
$ Z $.
By induction on stratum, we may assume that the assertion
holds for any
$ y \in \bC^{n}\setminus \{0\} $.
Note that the
$ b $-function of a global defining equation of a central
hyperplane arrangement is equal to that of a local equation at
$ 0 $, using the
$ \bC^{*} $-action.
We can apply 5.1 with
$ \alpha_{i} = 0 $ for any
$ i $, and (5.1.1) holds by [3] where
$ \omega\wedge = 0 $.
In particular,
$ H^{n-1}(\bP^{n-1}\setminus Z,\bC) $ is nonzero and is generated
by logarithmic forms on an embedded resolution of
$ (\bP^{n-1},Z) $, see [14].
Then
$ \Gr_{F}^{p}H^{n-1}(F_{0},\bC)_{1} = 0 $ for
$ p \ne n-1 $, and hence
$ -1 $ is the only integral root by Theorem~2.
Moreover,
$ \Gr^{W}_{i}H^{n-1}(F_{0},\bC)_{1} = 0 $ for
$ i \ne 2n-2 $ by the Hodge symmetry of
$ \Gr^{W}_{i}H^{n-1}(F_{0},\bC)_{1} $.
Thus
$ \Gr^{W}_{2n-2}H^{n-1}(F_{0},\bC)_{1} \ne 0 $.
(This also follows from [10] in the case
$ Z $ is a divisor with normal crossings.)
So the assertion holds from Lemma~3.6 together with [18], [25].

\medskip\noindent
{\bf 5.4.~$ b $-Function of a generic hyperplane arrangement.}
The
$ b $-function
$ b_{f}(s) $ of a generic central hyperplane arrangement with a 
reduced equation is determined by U.~Walther [43]:

\medskip\noindent
$ (W) $\,\,\,
The roots of
$ b_{f}(s) $ are
$ -j/d $ for
$ n \le j \le 2d - 2 $, and the multiplicity of a root
$ \alpha $ is
$ 1 $ for
$ \alpha \ne -1 $
and is
$ n $ for
$ \alpha = -1 $, assuming
$ d > n $.

\medskip
(The last assertion on the multiplicity of
$ -1 $ was not proved in loc.~cit.)
Here generic means that a central hyperplane arrangement has
normal crossings outside the origin.
In particular, the arrangement is not the pull-back of an
arrangement in a strictly lower dimensional vector space since
$ d > n $.
Using the arguments in this paper, we can give another proof
of his theorem as follows:

By Proposition~5.2 using [3], [14] in the normal crossing case,
we first get

\medskip\noindent
(5.4.1)\,\,\, The roots of
$ b_{f}(-s) $ is strictly smaller than
$ 2 $.

\medskip\noindent
The assertion on the integral roots follows from Proposition~5.3.
For the non-integral roots, the multiplicity is always
$ 1 $ by Proposition~4.6.
Moreover, Theorem~2 and (5.4.1) imply that
$ 1 + k/d $ is a root of
$ b_{f}(-s) $ for
$ 1 \le k \le d - 2 $, and is not a root for
$ k = d - 1 $.
Indeed, by Corollary~1,
$ m_{\alpha} $ for
$ \alpha = k/d $ is strictly smaller than
$ \binom{d-2}{n-1} $ if
$ k \le d - 2 $, and they coincide if
$ k = d - 1 $.
Since
$ \dim H^{n-1}(F_{0},\bQ)_{\lambda} = \binom{d-2}{n-1} $ for
$ \lambda = \exp(-2\pi ik/d) $ with
$ 1 \le k \le d - 1 $ (see [6], [28]), these imply
$$
\aligned
&F^{n-1}H^{n-1}(F_{0},\bC)_{\lambda} \ne H^{n-1}(F_{0},\bC)_{\lambda}
\quad\text{if}\quad 1 \le k \le d - 2,
\\
&F^{n-1}H^{n-1}(F_{0},\bC)_{\lambda} = H^{n-1}(F_{0},\bC)_{\lambda}
\quad\text{if}\quad k = d - 1.
\endaligned
$$
Note that
$ P^{n-1}H^{n-1}(F_{0},\bC) = F^{n-1}H^{n-1}(F_{0},\bC) $ because
$ P_{0}\cL^{(k/d)} = F_{0}\cL^{(k/d)} $ (since
$ Z $ is a divisor with normal crossings).
We have moreover
$ P^{n-2}H^{n-1}(F_{0},\bC)_{\lambda} =
H^{n-1}(F_{0},\bC)_{\lambda} $ by (5.4.1).
So the assertion follows from Theorem~2.

\medskip\noindent
{\bf 5.5.~Example of a nongeneric hyperplane arrangement.}
With the notation of 5.1, assume
$ n = 3 $,
$ d = 7 $, and
$ h = (x^{2}-y^{2})(x^{2}-1)(y^{2}-1) $ so that
$ f $ is as in 3.4 (ii).
Then
$ 5/7 $ is a root of
$ b_{f}(-s) $ (although it is not a jumping coefficient).

Indeed, let
$ Z' = \{x^{2}-y^{2} = 0\} $,
$ Z'' = \{(x^{2}-1)(y^{2}-1) = 0\} $.
Then we can apply the argument in 5.1 to the case where
$ \alpha_{i} = -5/7 $ if
$ Z_{i} \subset Z' $, and
$ \alpha_{i} = 2/7 $ if
$ Z_{i} \subset Z'' $.
In this case we have
$$
\dim \cA_{h,\alpha}^{1} = 6,\,\,\,
\dim \cA_{h,\alpha}^{2} = 9,\,\,\,
\dim H^{2}(\cA_{h,\alpha}^{\ssbull}) = \chi(U) = 4.
$$
Since
$ \Omega_{Y}^{2}\otimes_{\cO}P_{0}\cL^{(5/7)} \simeq \cO_{Y}(2) $
where
$ Y = \bP^{2} $, we see that
$ g(x,y)h^{-5/7}dx\wedge dy $ can be extended to a section of
$ \Omega_{Y}^{2}\otimes_{\cO}P_{0}\cL^{(5/7)} $ if
$ g(x,y) $ is a polynomial of degree
$ \le 2 $.
Moreover,
$ g(x,y)h^{-5/7}dx\wedge dy $ is contained in
$ \cA_{h,\alpha}^{2} $ if
$ g(x,y) $ is a linear combination of
$ (x - \varepsilon)(y - \varepsilon') $ with
$ \varepsilon, \varepsilon' = \pm 1 $, i.e. if
$ g(x,y) $ has order
$ \le 1 $ for both
$ x $ and
$ y $.
Indeed,
$ h^{\alpha} $ is naturally extended to a section of
$ P_{0}\cL^{(5/7)} $ having a simple zero along
$ Z'' $ and the divisor at infinity, and
$ \frac{d(x+\varepsilon)}{x+\varepsilon}\wedge
\frac{d(y+\varepsilon')}{y+\varepsilon'} $ has a simple pole
along the divisor at infinity.

Let
$ V $ be the vector subspace of
$ \cA_{h,\alpha}^{2} $ consisting of such elements.
We see that the dimension of the image of
$ d\cA_{h,\alpha}^{1} $ in
$ \cA_{h,\alpha}^{2}/V $ is at least
$ 2 $, calculating the differential
$ d $ of
$ \cA_{h,\alpha}^{\ssbull} $ which is defined by
$ \omega\wedge $, see 5.1.
So we get
$ \dim V\cap d\cA_{h,\alpha}^{1} < \dim V = 4 $ because
$ \dim d\cA_{h,\alpha}^{1} = 5 $.
Thus the image of
$ V $ in
$ H^{2}(U,L^{(5/7)}) $ does not vanish by (5.1.1) (using [14]),
and hence
$ 5/7 $ is a root of
$ b_{f}(-s) $ by Theorem~2 together with Proposition~4.9.

\medskip\noindent
{\bf 5.6.~Spectrum of a generic hyperplane arrangement.}
The calculation in 5.4 implies that the coefficient
$ m_{\alpha} $ of
$ \Sp(f,0) $ for
$ \alpha \in \bZ $ is given by
$$
m_{n-i} = (-1)^{n-1-i}
\hbox{$\binom{d-1}{i}$}\quad\text{for}\,\,
1 \le i \le n-1.
\leqno(5.6.1)
$$
Assume
$ \alpha \notin \bZ $.
Then
$ m_{\alpha} $ is calculated by Corollary~1 if
$ \alpha < 1 $ or
$ \alpha > n-1 $.
For
$ 1 < \alpha < n-1 $, it is possible to calculate
$ m_{\alpha} $ using [10] together with the (twisted)
weight spectral sequence
because the dimension of the cohomology of the twisted forms
$ \Omega_{\bP^{i}}^{j}(r) $ on the projective space
$ \bP^{i} $ can be calculated by
using the Bott vanishing theorem and the Euler sequence.

\medskip\noindent
{\bf 5.7.~Remark.}
For hyperplane arrangements, it is conjectured by Mustata [26]
that the jumping coefficients depend only on the combinatorial
data (i.e. the dimensions of various intersections of irreducible
components) of the hyperplane arrangement.
This assertion can be reduced to the one for the spectrum,
and will be proved in a forthcoming joint paper with Budur and
Mustata.
Using [5] together with Hodge theory, it is easy to show that
they remain unchanged under a deformation with the combinatorial
data fixed, see also [30].
However, the parameter space of hyperplane arrangements with fixed
combinatorial data is not connected as shown in [31].
In the case of a cone of a curve of higher degree in
$ \bP^{2} $,
a similar fact is known as Zariski's example, see e.g. [13].

For hyperplane arrangements, it is possible to show
the non-connectivity of the parameter space by using the following:

\medskip\noindent
(A)\,\, Let
$ p_{i} = a_{i} + \lambda b_{i}\,(i = 1, 2, 3) $ be three points on
$ \bC^{2} $ with a linear motion parametrized by
$ \lambda \in \bC $.
Then there are, in general, two values of
$ \lambda $ for which the three points are on one line.

\medskip
Indeed, this implies that, for a certain family of line
arrangements in
$ \bC^{2} $ with fixed combinatorial data whose parameter space is
one-dimensional, it is possible only for two points of the parameter
space to add one line to the corresponding line arrangement
so that the obtained line arrangement has certain restricted
combinatorial data.
For example, consider the union of
$$
\{xy(x^{2}-1)(y^{2}-1)(x-y)(x-y-1) = 0\} \subset \bC^{2},
$$
with three lines
$ L_{1}, L_{2}, L_{3} $ such that
$ L_{1} $ passes
$ (1,0) $,
$ L_{2} $ is parallel to
$ \{x = y\} $,
$ L_{3} $ passes
$ (0,0) $ and
$$
L_{1}\cap L_{2} \subset \{y = 1\},\quad
L_{2}\cap L_{3} \subset \{x = -1\},\quad
L_{1}\cap L_{3} \subset \{y = -1\}.
$$
The parameter space of such arrangements is one-dimensional if
$ L_{3} $ is deleted.
So we can apply the above argument to the three points
$ (0,0) $,
$ L_{1}\cap\{y = -1\} $, and
$ L_{2}\cap\{x = -1\} $.

\end{document}